\documentclass[a4paper,11pt]{amsart}  
\usepackage[all]{xy}
\usepackage{stmaryrd}
\usepackage{hyperref}
\usepackage{enumerate}
\usepackage{amssymb}
\newcommand{\G}{\mathcal{G}}
\newcommand{\s}{\mathbf{s}}
\renewcommand{\t}{\mathbf{t}}

\newcommand{\set}[1]{\left\{#1\right\}}

\newcommand{\eps}{\varepsilon}

\newcommand{\X}{\ensuremath{\mathfrak{X}}}
\newcommand{\F}{\ensuremath{\mathcal{F}}}

\newcommand{\al}{\alpha}
\newcommand{\be}{\beta} 

\newcommand{\Lie}{\mathcal{L}}  
\renewcommand{\gg}{\mathfrak{g}}
\newcommand{\hh}{\mathfrak{h}}

\newcommand{\tto}{\rightrightarrows}
\renewcommand{\Vert}{\text{\rm Vert}\,}

\renewcommand{\d}{\mathrm{d}}

\DeclareMathOperator{\Graph}{Graph}

\DeclareMathOperator{\Diff}{Diff}

\DeclareMathOperator{\Ker}{Ker}

\DeclareMathOperator{\Mor}{Mor}

\newtheorem{theorem}{Theorem}[section]

\newtheorem{proposition}{Proposition}[section]
\newtheorem{corollary}{Corollary}[section]

\theoremstyle{plain}
\newtheorem{definition}{Definition}[section]
\newtheorem{example}{Example} [section]

\newtheorem{remark}{Remark}[section]

\begin{document}

\title[The Symplectization Functor]{The Symplectization Functor}%

\author[R.L.~Fernandes]{Rui Loja Fernandes}

\address{Departamento de Matem\'{a}tica, 
Instituto Superior T\'{e}cnico, 1049-001 Lisboa, PORTUGAL} 

\begin{abstract}
We describe a symplectization functor from the Poisson category to the
symplectic ``category'' and we study some of its properties.
\end{abstract}


\maketitle

\section{Introduction}

Recall that to every Poisson manifold $(M,\pi)$ there is associated a Lie algebroid
structure on $T^*M$. It was shown in \cite{CaFe2,CrFe1,CrFe2} that we can always associate
to this Lie algebroid a topological groupoid $\Sigma(M)$, which in favorable circumstances
(the \emph{integrable case}) is a differentiable groupoid. In the integrable case, this Lie 
groupoid carries a compatible symplectic structure which, as we recall below, arises from the canonical 
symplectic structure on $T^*M$. In this paper we look into the 
functorial properties of this construction. Namely, we take the point of view that $\Sigma$ 
is a functor from the \emph{Poisson category} to the \emph{symplectic ``category''}. This may be 
useful, for example, in the problem of (geometric, deformation) quantization of
Poisson manifolds.

First of all, by the \emph{Poisson category} we mean the category \textbf{Poiss} 
whose objects are the (integrable) Poisson manifolds and whose morphisms are the Poisson maps. 
On the other hand, by the \emph{symplectic ``category''} we mean the ``category'' \textbf{Symp} 
whose objects are the symplectic manifolds and whose morphisms are the 
canonical relations. This ``category'' was introduced by Alan Weinstein 
in \cite{Wein1}, who also used the quotations marks as a reminder that this is not really a category, since 
composition of canonical relations is not always defined. This ``category'' contains a
\emph{symplectic groupoid ``subcategory''}, which we denote \textbf{SympGrp}. Now, as we will recall below, 
$\Sigma$ associates to a Poisson manifold $(M,\pi)$ a \emph{symplectic groupoid} $(\Sigma(M),\omega)$ and to a 
Poisson map $\phi:(M_1,\pi_1)\to (M_2,\pi_2)$ a canonical relation
(in fact, a \emph{Lagrangian subgroupoid}) $\Sigma(\phi)\subset \Sigma(M_1)\times\overline{\Sigma(M_2)}$, 
in such a way that composition of Poisson maps corresponds to composition of canonical relations. 
Therefore, $\Sigma$ is a functor from \textbf{Poiss} to \textbf{SympGrp}.

In the Poisson category one has many geometric constructions:
\begin{itemize}
	\item Passing to sub-objects;
	\item Passing to quotients;
	\item Forming fibrations;
	\item Forming connected sums;
	\item (...)
\end{itemize}
It is natural to look at what happens to each such geometric construction upon
applying the functor $\Sigma$. In this paper we survey recent
results and ongoing research on this direction.

The \emph{symplectization functor} $\Sigma$ should not be confused with the 
\emph{integration functor} $\G$ which goes from the category of Lie algebroids to the
category of Lie groupoids. For example, a Poisson map is almost never a Lie algebroid morphism
of the underlying cotangent algebroids. Symplectization of many of the 
constructions above entails a rich geometry, which is not present in Lie algebroid theory, and
which makes this passage from the Poisson category to the symplectic groupoid ``category'' 
far from being obvious.

The rest of this paper is organized into two sections. In the first section, we review
the construction of $\Sigma$. In the second section, we consider the effect
of $\Sigma$ on some of the geometric constructions mentioned above.
\vskip 10 pt

\begin{remark}
In order to simplify the presentation we will assume throughout this paper that our 
Poisson manifolds are integrable, However, as the reader will notice, many of the 
constructions below still make sense for non-integrable Poisson manifolds.
\end{remark}

\section{The $\Sigma$ functor}
Let $(M,\{\cdot,\cdot\})$ be a Poisson manifold. We will denote by
$\pi\in\X^2(M)$ the associated Poisson tensor which is given by
$\pi(\d f,\d  g):=\{f,g\}$, $(f,g \in  C^\infty(M))$,
and by $\pi^\sharp:T^\ast M \rightarrow TM$ the vector bundle map defined by
\[
\pi^\sharp (\d h)=X_h:=\{h,\cdot\}.
\]
As usual, we call $X_h$ the hamiltonian vector field determined by $h\in
C^\infty(M)$. Also, there is a Lie bracket on 1-forms 
$[~,~]:\Omega^1(M)\times\Omega^1(M)\to\Omega^1(M)$ which is defined by:
\[ [\al,\be]=\Lie_{\pi^\sharp\al}\be-\Lie_{\pi^\sharp\be}\al-\d\pi(\al,\be).\]
The triple $(p:T^*M\to M,[~,~],\pi^\sharp)$ is a Lie algebroid, called the 
\emph{cotangent Lie algebroid} of the Poisson manifold $(M,\pi)$.

\subsection{$\Sigma$ on objects}
Let us recall briefly the construction of the groupoid $\Sigma(M)$ canonically
associated with the Poisson manifold $(M,\pi)$ (more details can be found in \cite{CrFe2,CrFe1}):
\[ \Sigma(M):=\frac{\{\text{cotangent paths }\}}{\{\text{cotangent homotopies}\}}.\]
where:
\begin{itemize}
	\item A \emph{cotangent path} is a $C^1$-path $a:I\to T^*M$ such that 
	\[ \pi^\sharp(a(t))=\frac{\d}{\d t}p(a(t)). \]
	\item Two cotangent paths $a_0$ and $a_1$ are \emph{cotangent homotopic} if there
	exists a family of cotangent paths $a_\eps(t)=a(\eps,t)$, $\eps\in [0,1]$, such 
	that the solution $b=b(\epsilon,t)$ of the differential equation
    \begin{equation}\label{diffeq}
      \partial_{t}b-\partial_{\epsilon} a= T_{\nabla}(a, b),\quad
      b(\epsilon,0)=0,
    \end{equation}
 	satisfies $b(\eps,1)=0$ (\footnote{Here $T_\nabla$ denotes the torsion of a 
 	connection $\nabla$, while $\partial_t$ and $\partial_\epsilon$ are the covariant
 	derivatives along the $t$ and $\epsilon$ directions. One can show that this condition 
 	does not depend on the choice of $\nabla$.}).
\end{itemize}
\vskip 10 pt

There is a natural groupoid structure on $\Sigma(M)$:
\begin{itemize}
	\item[-] The \emph{source} and \emph{target} maps $\s,\t:\Sigma(M)\to M$ are given 
		by $\s([a])=p(a(0))$ and $\t([a])=p(a(1))$;
	\item[-] The \emph{multiplication} in $\Sigma(M)$ is defined by
	$[a_1]\cdot [a_0]=[a_1\odot a_0]$, where $a_1\odot a_0$ denotes concatenation of cotangent
	paths:
		\[
			a_{1}\odot a_{0}(t)\equiv \left\{
			\begin{array}{ll}
			2a_{0}(2t),\qquad& 0\le t\le \frac{1}{2},\\ \\
			2a_{1}(2t-1),\qquad & \frac{1}{2}< t\le 1.
			\end{array}
			\right.
		\]
	\item[-] The \emph{identity} section $\eps:M\to\Sigma(M)$ is given by
	$\eps(x)=[0_x]$, where $0_x$ is the trivial cotangent path based at $x$;
	\item[-] The \emph{inverse} map $\iota:\Sigma(M)\to\Sigma(M)$ is defined
	by $\iota([a])=[\overline{a}]$, where $\overline{a}(t)=a(1-t)$ denotes the
	opposite path.
\end{itemize}
The space of cotangent paths $P_\pi(M)$ is furnished with the $C^2$ compact-open 
topology, and in $\Sigma(M)=P_\pi(M)/\sim$ we take the quotient topology. 
Then $\Sigma(M)$ becames a topological groupoid, which is sometimes called the \emph{Weinstein
groupoid} of $(M,\pi)$. In general, this groupoid is not smooth, but it is rather 
a differentiable stack (see \cite{TsZh1,TsZh2}). However, if $\Sigma(M)$ is smooth 
then one obtains a source 1-connected Lie groupoid that integrates the 
cotangent Lie algebroid of $(M,\pi)$, and in this case we
call $(M,\pi)$ an \emph{integrable Poisson manifold}. The precise obstructions
to integrability were determined in \cite{CrFe2,CrFe1}.

Let us now recall:

\begin{definition}
A symplectic form $\omega$ in a groupoid $\G$ is
\emph{multiplicative} iff the graph of the multiplication 
\[ \gamma_m=\{(g,h,g\dot h)\in \G\times\G\times\G~|~(g,h)\in\G^{(2)}\}\] 
is a Lagrangian submanifold of $\G\times\G\times\bar{\G}$. The pair
$(\G,\omega)$ is called a \emph{symplectic groupoid}.
\end{definition}

If $(\G,\omega)$ is a symplectic groupoid then the base manifold $M$
has a canonical Poisson bracket such that:
\begin{enumerate}[(i)]
	\item $\s$ is Poisson and $\t$ is anti-Poisson;
	\item the Lie algebroid of $\Sigma$ is canonically isomorphic to $T^*M$.
\end{enumerate}
Conversely, we have:

\begin{theorem}[\cite{CaFe2,CrFe2}]
Let $(M,\pi)$ be an integrable Poisson manifold. Then $\Sigma(M)$ is 
a symplectic groupoid whose Poisson structure on the base coincides with $\pi$.
\end{theorem}

For the sequel, it is important to understand how the symplectic structure on
$\Sigma(M)$ arises. For that, we need to look closer at cotangent homotopies
and at an alternative description of these homotopies in terms of a 
Lie algebra action. 

We will denote by $P(T^*M)$ the space of all $C^1$-paths 
$a:I\to T^*M$. This is a Banach manifold 
in the obvious way, and the space of cotangent paths 
\[ P_\pi(M)\subset P(T^*M)\] 
is a Banach submanifold. Also, we let $P(M)$ denote the space
of $C^2$-paths $\gamma:I\to M$. A basic fact, which in fact explains the
existence of a symplectic structure on $\Sigma(M)$, is that 
\[ P(T^*M)=T^*P(M), \]
so that $P(T^*M)$ carries a natural (weak) symplectic structure
$\omega_{\text{can}}$. 

Let us identify the tangent space $T_{a}P(T^*M)$ with the space of vector fields
along $a$:
\[ T_{a}P(T^*M)=\set{U:I\to TT^*M ~|~ U(t)\in T_{a(t)}T^*M}.\]
and denote by $\omega_0$ the canonical symplectic form on $T^*M$. Then
the 2-form $\omega_{\text{can}}$ on $P(T^*M)=T^*P(M)$ is given by:
\begin{equation}
  \label{eq:lifted:symp:form}
  (\omega_{\text{can}})_a(U_1, U_2)=\int_0^1 \omega_0(U_1(t), U_2(t))dt,
\end{equation}
for all $U_1, U_2\in T_a P(T^*M)$. Moreover, it is easy to check that
$\d\omega_{\text{can}}=0$ and that
$\omega_{\text{can}}^\sharp:T_{a}P(T^*M)\to T_{a}^*P(T^*M)$ is 
injective. Hence, $\omega_{\text{can}}$ is a weak symplectic form. 

On the other hand, the Lie algebra
\[
P_{0}\Omega^1(M):=\set{\eta_{t}\in \Omega^1(M), t\in I~|~\eta_{0}
                =\eta_{1}=0,\ \eta_t\text{ of class }C^1\text{ in }t}
\]
with the pointwise Lie bracket, acts on $P(T^*M)$ in such a way that:
\begin{enumerate}[(a)]
\item the action is tangent to $P_\pi(M)$;
\item two cotangent paths are homotopic if and only if they belong to
the same orbit.
\end{enumerate}
Now we have the following remarkable fact, first observed in \cite{CaFe2}:

\begin{theorem}
  \label{thm:sigma:model}
  The infinitesimal action of $P_{0}\Omega^1(M)$ on
  $(P(T^*M),\omega_{\text{can}})$ is Hamiltonian, with equivariant
  moment map $J:P(T^*M)\to P_{0}\Omega^1(M)^*$ given by
  \begin{equation}
  \label{eq:moment:map}
  \langle J(a),\eta\rangle=
  \int_0^1\langle\frac{d}{dt}\pi(a(t))-\#a(t),\eta(t,\gamma(t))\rangle\,dt.
\end{equation}
\end{theorem}

Therefore, the groupoid $\Sigma(M)$ is obtained by symplectic reduction:
\[ \Sigma(M)=P(T^*M)//P_0\Omega^1(M)=J^{-1}(0)/P_0\Omega^1(M).\]
This gives the multiplicative symplectic form on $\Sigma(M)$.

\begin{remark}
\label{rem:symp}
It is important to note that not every Lie groupoid integrating the cotangent bundle
$T^*M$ of a Poisson manifold is a symplectic groupoid. What we saw above is
that the unique source 1-connected groupoid is indeed a symplectic groupoid, and that 
this symplectic structure arises from the canonical symplectic structure on $T^*M$.
\end{remark}

\subsection{$\Sigma$ on morphisms}
Now that we know what the effect of $\Sigma$ on objects is, let us look
at its effect on a Poisson morphism 
$\phi:(M_1,\pi_1)\to (M_2,\pi_2)$. Note that, in general, $\phi$ \emph{does not} induce
a morphism of Lie algebroids, so the answer \emph{is not} a Lie groupoid morphism 
$\Sigma(\phi):\Sigma(M_1)\to\Sigma(M_2)$. In order to find out what the answer should be, let us
recall different ways of expressing the condition for a map to be Poisson:

\begin{proposition}[\cite{Wein2}]
Let $\phi:(M_1,\pi_1)\to (M_2,\pi_2)$ be a smooth map between two Poisson manifolds. 
The following conditions are equivalent:
\begin{enumerate}[(a)]
	\item The map $\phi$ preserves Poisson brackets: $\{f\circ\phi,g\circ\phi\}_1=\{f,g\}_2\circ\phi$.
	\item The Poisson bivectors are $\phi$-related: $\phi_*\pi_1=\pi_2$.
	\item $\Graph(\phi)\subset M_1\times\overline{M}_2$ is a coisotropic submanifold\emph{
	(\footnote{A submanifold $N$ of a Poisson manifold $(M,\pi)$ is called \emph{coisotropic} if
	$\pi^\sharp(TN)^0\subset TN$, where $(TN)^0\subset T^*M$ denotes the annihilator of $TN$. Also,
	$\overline{M}$ denotes the Poisson manifold $M$ with the symmetric Poisson structure.})}.
\end{enumerate}
\end{proposition}

It is clear that if $(M_1,\pi_1)$ and $(M_2,\pi_2)$ are Poisson manifolds, then we have 
\[ \Sigma(M_1\times\overline{M}_2)=\Sigma(M_1)\times\overline{\Sigma(M_2)}.\]
Hence, to integrate a Poisson morphism, we just need to know what objects integrate
coisotropic submanifolds of a Poisson manifold. This problem was solved by Cattaneo and
Felder in \cite{CaFe1}:

\begin{theorem}
If $\G\tto C$ is a Lagrangian subgroupoid of a symplectic groupoid $\Sigma\tto M$, 
then $C\subset M$ is a coisotropic submanifold. Conversely, if $C$ is a coisotropic
submanifold of an integrable Poisson manifold $(M,\pi)$, then there exists a Lagrangian
subgroupoid $\G\tto C$ of $\Sigma(M)$ that integrates $C$. 
\end{theorem}

Let us explain why coisotropic submanifolds integrate to Lagrangian
subgroupoids. Note that $C$ is a coisotropic submanifold of $(M,\pi)$ iff
its conormal bundle $\nu^*(C):=(TC)^0\subset T^*M$ is a Lie subalgebroid
of the cotangent Lie algebroid $T^*M$. Therefore, if $(M,\pi)$ is 
integrable, then there exists a source connected Lie subgroupoid $\G\tto C$
of the groupoid $\Sigma(M)\tto M$ that integrates $\nu^*(C)$. Now we claim
that the restriction of the symplectic form $\omega$ to $\G$ vanishes which,
combined with $\dim\G=1/2\dim\Sigma(M)$, implies that $\G$ is Lagrangian.

To prove our claim, we observe that for the canonical symplectic form $\omega_0$ 
on $T^*M$ the submanifold $\nu^*(C)\subset T^*M$ is Lagrangian. Using the explicit
expression (\ref{eq:lifted:symp:form}) for the symplectic form $\omega_{\text{can}}$
on $P(T^*M)$, we see immediately that space of paths $P(\nu^*(C))\subset P(T^*M)$ is 
isotropic. It follows that the symplectic form $\omega$ on the symplectic quotient 
$\Sigma(M)$ restricts to zero on the submanifold $\G=P(\nu^*(C))//P_0\Omega^1(M)$, 
as we claimed.

We have now found what Poisson morphisms integrate to:

\begin{corollary}
\label{cor:mor}
Let $\phi:(M_1,\pi_1)\to (M_2,\pi_2)$ be a Poisson map between two integrable
Poisson manifolds. Then $\phi$ integrates to a Lagrangian subgroupoid
$\Sigma(\phi)\subset \Sigma(M_1)\times\overline{\Sigma(M_2)}$.
\end{corollary}

\subsection{The symplectic ``category''}
Let us recall now the symplectic ``category'' of Alan Weinstein \cite{Wein1}, 
which we will denote by \textbf{Symp}.

In the category \textbf{Symp} the objects are the symplectic manifolds and
the morphisms are the \emph{canonical relations}. We recall that if 
$(S_1,\omega_1)$ and $(S_2,\omega_2)$ are two symplectic manifolds, 
then a \emph{canonical relation} from $S_1$ to $S_2$ is, by definition, 
a Lagrangian submanifold $L\subset S_1\times\overline{S_2}$ where, as usual, 
we denote by a bar the same manifold with the opposite symplectic 
structure.

There is, however, a problem: if $L_1\in \Mor(S_1,S_2)$ and $L_2\in\Mor(S_2,S_3)$ 
are canonical relations, their composition:
\[ L_1\circ L_2:=\{(x,z)\in S_1\times S_3~|~\exists y\in S_2,\text{ with }(x,y)\in L_1
\text{ and }(y,z)\in L_2\},\]
may not be a smooth submanifold of $S_1\times \overline{S_3}$. One needs a certain 
clean intersection property which we will not discuss here (see \cite{Wein1}).
However, whenever it is a smooth submanifold, it is indeed a Lagrangian submanifold.
Therefore, in \textbf{Symp} composition is not always defined and Weinstein proposed to
name it a ``category'', with quotation marks. 

It maybe worth to point out three simple properties of the symplectic ``category''. First, 
in the symplectic ``category'' the points of an object $(S,\omega)$ are the elements 
of $\Mor(\text{pt},S)$, i.e., the Lagrangian submanifolds $L\subset S$. Second, there
exists an involution in \textbf{Symp}: it is the covariant functor of \textbf{Symp}
which takes an object $S$ to itself $S^\dag:=S$ and a morphism $L\in\Mor(S_1,S_2)$ to the 
morphism $L^\dag:=\{(y,x)\in S_2\times\overline{S_1}~|~(x,y)\in L\}\in\Mor(S_2,S_1)$.
Finally, in \textbf{Symp} we have the ``subcategory'' of symplectic groupoids, denoted 
\textbf{SympGrp}, where the objects are the symplectic groupoids $(\G,\omega)$ and the
morphisms $\Mor(\G_1,\G_2)$ are the Lagrangian subgroupoids of $\G_1\times\overline{\G_2}$. 

Now let us denote by \textbf{Poiss} the Poisson category, in which the objects are
the Poisson morphisms and the morphisms are the Poisson maps. We can summarized the two 
previous paragraphs above by saying that $\Sigma$ is a covariant functor from \textbf{Poiss} to 
\textbf{SympGrp}. 

\subsection{Examples} For illustration purposes we give three, well-known, classes of examples.

\subsubsection{Symplectic manifolds.} If $(M,\pi)$ is a non-degenerate Poisson structure, then
$\omega=(\pi)^{-1}$ is a symplectic form. In this case, a cotangent path $a:I\to T^*M$
is completely determined by its base path $\gamma:I\to M$, since we have:
\[ a(t)=\omega(\dot{\gamma}(t),\cdot).\]
Moreover, two cotangent paths $a_0$ and $a_1$ are cotangent homotopic iff their base paths
$\gamma_0$ and $\gamma_1$ are homotopic relative to its end-points. Hence, $\Sigma(M)$
is just the fundamental groupoid of $M$, and so can be identified as:
\[ \Sigma(M)=\widetilde{M}\times_{\pi_1(M)}\widetilde{M},\]
where $\widetilde{M}$ is the universal covering space and $\pi_1(M)$ acts diagonally. This 
space carries an obvious symplectic groupoid structure which is inherit from the pair 
(symplectic) groupoid $\widetilde{M}\times\widetilde{M}$. Note that when $(M_1,\pi_1)$ 
and $(M_2,\pi_2)$ are both symplectic, a Poisson map $\phi:M_1\to M_2$ is not, in general, 
symplectic (such a Poisson map $\phi$ must be a submersion).

\subsubsection{Trivial Poisson manifolds.} A smooth manifold $M$ always carries the 
trivial Poisson structure $\pi\equiv 0$.  In this case, a cotangent path is just a 
path into a fiber of $T^*M$, and any such path $a:I\to T_x^*M$ is cotangent homotopic
to its average $\bar{a}(t):=\int_0^1 a(s)\d s$ (a constant path in $T^*_xM$). It is easy
to check that two cotangent paths are homotopic iff their averages are equal, so that
$\Sigma(M)=T^*M$. On the other hand, any smooth map $\phi:M_1\to M_2$ is a 
Poisson map and we have $\Sigma(\phi)=\text{Graph}(\phi^*)\subset T^*M_1\times \overline{T^*M_2}$.
It follows that $\Sigma$ embeds the category of smooth manifolds in \textbf{Symp}.

\subsubsection{Linear Poisson structures.} It is well known (according to \cite{Wein}, it
was already known to Lie!) that a linear Poisson structure on a vector space is the same thing
as a Lie algebra structure on the dual vector space. If $\gg$ is a Lie algebra, we let 
$M=\gg^*$ with is canonical linear Poisson structure. Then $\Sigma(\gg^*)=T^*G$, where $G$ 
is the 1-connected Lie group with Lie algebra $\gg$. Here the symplectic structure is 
the canonical symplectic structure on the cotangent bundle $T^*G$ and the groupoid 
structure can be defined as follows: let $\s,\t:T^*G\to\gg^*$ be the left/right trivializations: 
\[ \s(\al_g)=(\d_eL_g)^*\al_g,\quad \t(\be_h)=(\d_eR_h)^*\be_h,\quad (\al_g\in T_g^*G, 
\be_h\in T_h^*G).\]
Then, if $\al_g,\be_h\in T^*G$ are composable (i.e., $\s(\al_g)=\t(\be_h)$) then their product is 
given by:
\[ \al_g\cdot\be_h:=(d_e L_g)^*\be_h\in T_{gh}^*G.\]
If $\phi:\gg\to\hh$ is a Lie algebra homomorphism, then $\phi^*:\hh^*\to\gg^*$ is a 
Poisson map. Upon applying the symplectization functor we obtain a Lagrangian submanifold
$\Sigma(\phi^*)\subset T^*H\times \overline{T^*G}$. It is not hard to check that if $\Phi:G\to H$
is the Lie group homomorphism with $\d_e\Phi=\phi$, then:
\[ \Sigma(\phi^*)=\text{graph}(\Phi^*)=\{(\be_{\Phi(g)},\al_g)\in T^*H\times \overline{T^*G}~|~(\d_g\Phi)^*\be_{\Phi(g)}=\al_g\}.\]
In this way, $\Sigma$ embeds the category of Lie algebras \textbf{LieAlg} in the 
category \textbf{SympGr}.

\section{Symplectization of Poisson geometry}
Now that we know how the functor $\Sigma$ is constructed, we can look at is
effect on various geometric constructions in Poisson geometry.

\subsection{Sub-objects}
A sub-object in the Poisson category is just a Poisson submanifold $N$ of 
a Poisson manifold $(M,\pi)$ (\footnote{Actually, this notion of sub-object
is too restrictive, and in \cite{CrFe1} it is explained how the class of, so called, 
\emph{Poisson-Dirac submanifolds} is the more appropriate class of sub-objects to consider.}). 
However, in general, a Poisson submanifold of an integrable Poisson manifold is not integrable (see 
\cite{CrFe2,Xu1}). We will assume that both $(M,\pi)$ and $N$ are integrable.
Then the inclusion $i:N\hookrightarrow M$ is a Poisson morphism which, according
to Corollary \ref{cor:mor}, integrates to a Lagrangian subgroupoid $\Sigma(i)\subset\Sigma(N)\times\overline{\Sigma(M)}$. However, this is not
the end of the story.

For a Poisson submanifold $N\subset M$, let us consider the set of equivalence classes of
cotangent paths that take their values in the restricted subbundle $T^*_NM$:
\[ \Sigma_N(M):=\{[a]\in\Sigma(M)~|~a:I\to T^*_NM\}.\]
This is a Lie subgroupoid of $\Sigma(M)$: it is the Lie subgroupoid that integrates 
the Lie subalgebroid $T^*_NM\subset T^*M$. Moreover, this subalgebroid is a coisotropic 
submanifold of the symplectic manifold $T^*M$, and it follows that $\Sigma_N(M)\subset \Sigma(M)$
is a coisotropic Lie subgroupoid. The fact that the closed 2-form $\omega$ on $\Sigma_N(M)$ is multiplicative
implies that if we factor by its kernel foliation, we still obtain a symplectic groupoid, and in fact
(see \cite{CrFe1}):
\[ \Sigma(N)\simeq \Sigma_N(M)/\Ker\omega.\]

These two constructions are related as follows:

\begin{theorem}
Let $i:N\hookrightarrow M$ be an integrable Poisson submanifold of an 
integrable Poisson manifold, and $\Sigma(i)\subset\Sigma(N)\times\overline{\Sigma(M)}$
the corresponding Lagrangian subgroupoid. For the restriction of the projections on each factor:
\[
\xymatrix{
&\Sigma(i)\ar[dl]_{\pi_1}\ar[dr]^{\pi_2}\\
\Sigma(N)&&\overline{\Sigma(M)}}\]
$\pi_2$ is a diffeomorphism onto the coisotropic subgroupoid
$\Sigma_N(M)\subset \Sigma(M)$ above, and the groupoid morphism:
\[ (\pi_1)\circ(\pi_2)^{-1}:\Sigma_N(M)\to \Sigma(N),\]
corresponds to the quotient map $\Sigma_N(M)\to\Sigma_N(M)/\Ker\omega\simeq\Sigma(N)$.
\end{theorem}

A very special situation happens when the exact sequence of Lie algebroids
\[ \xymatrix{0\ar[r]& (TN)^0\ar[r]& T^*_NM\ar[r]& T^*N\ar[r]& 0}\]
splits: in this case, the splitting $\phi:T^*M\to T^*_NM\subset T^*M$ 
integrates to a symplectic Lie groupoid homomorphism 
$\Phi:\Sigma(N)\to\Sigma(M)$, which realizes $\Sigma(N)$ as a symplectic subgroupoid
of $\Sigma(M)$. Poisson submanifolds of this sort maybe called \emph{Poisson-Lie submanifolds},
and there are topological obstructions on a Poisson submanifold for this to happen (see
\cite{CrFe1}).

\subsection{Quotients}

What we have just seen for sub-objects is typical: though the functor $\Sigma$ gives us
some indication of what the integration of a certain geometric construction is, there is
often extra geometry hidden in the symplectization. Another instance of this happens when 
one looks at quotients.

Let $G$ be a Lie group that acts smoothly by Poisson diffeomorphisms on a Poisson manifold 
$(M,\pi)$. We will denote the action by $\Psi:G\times M\to M$ and we will also write 
$\Psi(g,x)=g\cdot x$. For each $g\in G$, we set:
\[ \Psi_g:M\to M,\ x\mapsto g\cdot x,\]
so that each $\Psi_g$ is a Poisson diffeomorphism. 

Now we apply the functor $\Sigma$. For each $g\in G$, we obtain a Lagrangian
subgroupoid $\Sigma(\Psi_g)\subset \Sigma(M)\times\overline{\Sigma(M)}$. This 
Lagrangian subgroupoid is, in fact, the graph of a symplectic Lie groupoid
automorphism, which we denote by the same symbol $\Sigma(\Psi_g):\Sigma(M)\to\Sigma(M)$.
Also, it is not hard to check that 
\[ \Sigma(\Psi):G\times \Sigma(M)\to\Sigma(M), (g,[a])\mapsto g\cdot [a]:=\Sigma(\Psi_g)([a]),\]
defines a symplectic smooth action of $G$ on $\Sigma(M)$. Briefly, $\Sigma$
lifts a Poisson action $\Psi:G\times M\to M$ to a symplectic action 
$\Sigma(\Psi):G\times \Sigma(M)\to\Sigma(M)$ by groupoid automorphisms. However, this is not
the end of the story.

Let us look closer at how one lifts the action from $M$ to $\Sigma(M)$. First of all, recall 
that any smooth action $G\times M\to M$ has a lifted cotangent action $G\times T^*M\to T^*M$. 
This yields, by composition, an action of $G$ on cotangent paths: if $a:I\to T^*M$ is a 
cotangent path we just move it around
\[ (g\cdot a)(t):=g\cdot a(t).\]
The fact that the original action $G\times M\to M$ is Poisson yields that
(i) $g\cdot a$ is a cotangent path whenever $a$ is a cotangent path, and (ii) if
$a_0$ and $a_1$ are cotangent homotopic then so are the translated paths $g\cdot a_0$ and $g\cdot a_1$.
Therefore, we have a well-defined action of $G$ on cotangent homotopy classes and this is just the
lifted action: $\Sigma(\Psi_g)([a])=[g\cdot a]$.

Now we invoke a very simple (but important) fact from symplectic geometry: for any 
action $G\times M\to M$ the lifted cotangent action $G\times T^*M\to T^*M$ is 
a hamiltonian action with equivariant momentum map $j:T^*M\to\gg^*$ given by:
\[ j:T^*M\to\gg^*,\ \langle j(\al_x),\xi\rangle:=\langle \al_x,X_\xi(x)\rangle,\]
where $X_\xi\in\X(M)$ is the infinitesimal generator associated with $\xi\in\gg$.
This leads immediately to the fact that the lifted action $\Sigma(\Psi):G\times \Sigma(M)\to\Sigma(M)$
is also hamiltonian (\footnote{Recall again that, after all, the symplectic structure on
$\Sigma(M)$ comes from the canonical symplectic structure on $T^*M$.}). The 
equivariant momentum map $J:\Sigma(M)\to\gg^*$ for the lifted action is given by:
\begin{equation}
\label{eq:lifted:momentum:map}
\langle J([a]),\xi\rangle:=\int_0^1 j(a(t))\d t=\int_a X_\xi.
\end{equation}
Since each $X_\xi$ is a Poisson vector field, the last expression shows that only
the cotangent homotopy class of $a$ matters, and $J$ is indeed well-defined. 
Expression \ref{eq:lifted:momentum:map} mean that we can see the momentum map of 
the $\Sigma(\Psi)$-action in two ways:
\begin{itemize}
\item It is the integration of the momentum map of the lifted cotangent action;
\item It is the integration of the infinitesimal generators along cotangent paths.
\end{itemize}
In any case, expression (\ref{eq:lifted:momentum:map}) for the momentum map 
shows that it satisfies the following additive property:
\[ J([a_0]\cdot [a_1])=J([a_0])+J([a_1]).\]
Hence $J$ is a groupoid homomorphism from $\Sigma(M)$ to the additive group $(\gg^*,+)$
or, which is the same, $J$ is differentiable \emph{groupoid 1-cocycle}. Moreover, 
this cocycle is \emph{exact} iff there exists a map $\mu:M\to\gg^*$ such that
\[ J=\mu\circ\t-\mu\circ\s,\]
and this happens precisely iff the original Poisson action $\Psi:G\times M\to M$ is hamiltonian
with equivariant momentum map $\mu:M\to\gg^*$. These facts, in one form or another, can be found 
in \cite{CoDaMo,CoDaWe,MiWe,WeXu}. We summarize them:

\begin{theorem}
\label{thm:lifted:action}
Let $\Psi:G\times M\to M$ be a smooth action of a Lie group $G$ on a
Poisson manifold $M$ by Poisson diffeomorphisms. There exists a
lifted action $\Sigma(\Psi):G\time\Sigma(M)\to\Sigma(M)$ by symplectic groupoid
automorphisms. This lifted $G$-action is Hamiltonian and admits the
momentum map $J:\Sigma(M)\to\gg^*$ given by (\ref{eq:lifted:momentum:map}). 
Furthermore:
\begin{enumerate}[(i)]
\item The momentum map $J$ is $G$-equivariant and is a groupoid 1-cocycle.
\item The $G$-action on $M$ is hamiltonian with momentum map $\mu:M\to\gg^*$
if and only if $J$ is an exact cocycle.
\end{enumerate}
\end{theorem}

Let us now assume that the Poisson action $\Psi:G\times M\to M$ is 
proper and free. These assumptions guarantee that $M/G$ is a smooth manifold. The
space $C^\infty(M/G)$ of smooth functions on the quotient is naturally identified 
with the space $C^\infty(M)^G$ of $G$-invariant functions on $M$. Since the Poisson 
bracket of $G$-invariant functions is a $G$-invariant function, we have a quotient 
Poisson structure on $M/G$ such that the natural projection $M\to M/G$ is a Poisson map.
It can be shown, using the results of \cite{CrFe2}, that if $M$ is an integrable Poisson
structure then $M/G$ is also integrable, and so the question arises: what is the relationship
between the symplectic groupoids $\Sigma(M)$ and $\Sigma(M/G)$?

First notice that if the original Poisson action $\Psi:G\times M\to M$ is 
proper and free, so is the lifted action $\Sigma(\Psi):G\times\Sigma(M)\to\Sigma(M)$. 
Therefore $0\in\gg^*$ is a regular value of the momentum map $J:\Sigma(M)\to\gg^*$.
Let us look at the symplectic quotient:
\[ \Sigma(M)//G:=J^{-1}(0)/G. \]
Since $J$ is a groupoid homomorphism, its kernel $J^{-1}(0)\subset\Sigma(M)$
is a Lie subgroupoid. Since $J$ is $G$-equivariant, the action leaves $J^{-1}(0)$ invariant and 
the restricted action is a free action by groupoid automorphisms. Hence, the groupoid structure
descends to a groupoid structure $\Sigma(M)//G\tto M/G$. It is easy to check then that this is
indeed a symplectic groupoid. In general, however, it is not true that:
\[ \Sigma(M/G)=\Sigma(M)//G,\]
so, in general, symplectization \emph{does not} commute with reduction. First of all, 
$J^{-1}(0)$ may not be connected, so that $\Sigma(M)//G$ may not have source connected 
fibers. Even if we restrict to $J^{-1}(0)^c$, the connected component of the identity section 
(so that the source fibers of $\Sigma(M)//G$ are connected) these fibers may have a non-trivial
fundamental group. This problem, as well as other issues such as non-free actions, convexity, etc,
is the subject of ongoing research (see \cite{FeOrRa}).

\subsection{Fibrations}

Let $(F,\pi)$ be a Poisson manifold. We denote by $\Diff_{\pi}(F)$ the
group of Poisson diffeomorphisms of $F$. We are interested in the 
following class of fibrations:

\begin{definition}
A \emph{Poisson fibration} $p:M\to B$ is a locally trivial fiber
bundle, with fiber type a Poisson manifold $(F,\pi)$ and with
structure group a subgroup $G\subset \Diff_{\pi}(F)$. When $\pi$ is
symplectic the fibration is called a \emph{symplectic fibration}.
\end{definition}

If $p:M\to B$ is a Poisson fibration modeled on a Poisson manifold
$(F,\pi)$, each fiber $F_b$ carries a natural Poisson structure
$\pi_b$: if $\phi_i:p^{-1}(U_i)\to U_i\times F$ is a local
trivialization, $\pi_b$ is defined by:
\[ \pi_b=(\phi_i(b)^{-1})_*\pi,\]
for $b\in U_i$. It follows from the definition that this 2-vector
field is independent of the choice of trivialization. Note that the Poisson
structures $\pi_b$ on the fibers can be glued to a Poisson structure
$\pi_V$ on the total space of the fibration:
\[ \pi_V(x)=\pi_{p(x)}(x), \quad (x\in M).\]
This 2-vector field is \emph{vertical}: $\pi_V$ takes values in
$\wedge^2\Vert\subset\wedge^2TM$. In this way, the fibers $(F_b,\pi_b)$ 
became Poisson submanifolds of $(M,\pi_V)$.

\begin{example}
\label{ex:neighb:leaf}
An important class of Poisson fibrations is obtained as follows. Take
any Poisson manifold $(P,\pi)$, fix a symplectic leaf $B$ of $P$ which
is an embedded submanifold, and let $p:M\to B$ be a tubular
neighborhood of $B$ in $P$. Each fiber carries a natural Poisson
structure, namely, the transverse Poisson structure
(\cite{Wein}). These transverse Poisson structures are all Poisson
diffeomorphic and it follows from the Weinstein splitting theorem that
this is a Poisson fibration.
\end{example}

The usual method to construct fibrations from a principal $G$-bundle
and a $G$-manifold $F$, works also for Poisson fibrations. One 
specifies the following data:
\begin{enumerate}[(a)]
\item A Lie group $G$;
\item A principal (right) $G$-bundle $P\to B$;
\item A Poisson (left) action $G\times F\to F$.
\end{enumerate}
On $P$ we consider the zero Poisson structure, so that
$P\times F$ becomes a Poisson $G$-manifold for the product 
Poisson structure and the diagonal $G$-action. It follows that 
the associated fibration $M=P\times_G F\to B$ carries a 
(vertical) Poisson structure which makes it into a Poisson fibration.
 
If we allow infinite dimensional Lie groups, the converse
is also true: let $p:M\to B$ be any Poisson fibration
with fiber type a Poisson manifold $(F,\pi)$. We let $G=\Diff(F,\pi)$ 
be the group of all Poisson diffeomorphisms, and we define a
principal $G$-bundle, called the \emph{Poisson frame bundle}:
\[ P=\bigcup_{x\in b}P_b\to B\]
where each fiber $P_b$ is the set of all Poisson diffeomorphisms 
$u:F\to F_b$. The group $G$ acts on (the right of) $P$ by
pre-composition: 
\[ P\times G\to P:~(u,g)\mapsto u \circ g,\]
and our original Poisson fiber bundle is canonically isomorphic to
the associated fiber bundle: $M\simeq P\times_G F$.

In the sequel, we will always assume that we have a Poisson fibration 
$M=P\times_G F\to B$ with fiber type $(F,\pi)$. We now apply our 
functor $\Sigma$ to this fibration: assuming that $(F,\pi)$ is an 
integrable Poisson manifold, we have a symplectic groupoid $\Sigma(F)$ and
the Poisson action $G\times F\to F$ lifts to a hamiltonian action 
$G\times \Sigma(F)\to \Sigma(F)$, which is by groupoid automorphisms. 
Let us denote by $\Sigma(M):=P\times_G\Sigma(F)\to B$ the associated fiber
bundle. There are two things to note about this bundle:
\begin{enumerate}[(i)]
\item $\Sigma(M)\to B$ is a symplectic fibration: this follows from the method
explained above to construct Poisson fibrations, except that now our fiber type
is symplectic;
\item $\Sigma(M)\tto M$ is a groupoid: all the structure maps, as well as
the composition, are defined in the obvious way from the groupoid structure
on the fiber.
\end{enumerate}

The groupoid $\Sigma(M)\tto M$ is an example of a \emph{fibered groupoid}. 
Formally, this means an internal category in the category of
fibrations where every morphisms is an isomorphism. In practice, both 
the total space $\G$ and the base $M$ of a fibered groupoid 
are fibrations over $B$ and all structure maps are fibered maps. For example, 
the source and target maps are fiber preserving maps over the identity:
\[
\xymatrix{\G\ar@<.5ex>[rr]\ar@<-.5ex>[rr]\ar[dr]& &M\ar[dl]\\ & B}
\]
In particular, each fiber of $\G\to B$ is a groupoid over a fiber of
$M\to B$. Moreover, the orbits of $\G$ lie inside the fibers of the 
base $M$. 

The construction above is just an instance of a general procedure to 
construct fibered groupoids: Let $P\to B$ be a principle $G$-bundle 
and assume that $G$ acts on a groupoid $\F\tto F$ by groupoid automorphisms. 
Then the associated fiber bundles $\G=P\times_G\F$ and $M=P\times_G F$ are the
spaces of arrows and objects of a fibered groupoid. Clearly, every
fibered groupoid is of this form provided we allow infinite
dimensional structure groups. We will say that the fibered groupoid
$\G$ has \emph{fiber type} the groupoid $\F$. Now we set:

\begin{definition}
A \emph{fibered symplectic groupoid} is a fibered groupoid $\Sigma$
whose fiber type is a symplectic groupoid $(\F,\omega)$. 
\end{definition}

Therefore, if $\G$ is a fibered symplectic groupoid over $B$, then 
$\G\to B$ is a symplectic fibration, and each symplectic fiber $\F_b$ 
is in fact a symplectic groupoid over the corresponding fiber $F_b$ of
$M\to B$. Since the base of a symplectic groupoid has a natural
Poisson structure for which the source (respectively, the target) is a
Poisson (respectively, anti-Poisson) map, we can summarize our results
as follows:

\begin{theorem}
The base $M\to B$ of a fibered symplectic groupoid $\Sigma\tto M$ has a natural
structure of a Poisson fibration.  Conversely, for any Poisson 
fibration $M\to B$, with fiber type $(F,\pi)$ 
an integrable Poisson manifold, there exists a unique (up to isomorphism)
source 1-connected fibered symplectic groupoid $\Sigma(M)\tto M$ 
integrating $M\to B$. 
\end{theorem}

In a nut shell, $\Sigma$ takes Poisson fibrations to symplectic groupoid
fibrations, i.e., fibrations in \textbf{Poiss} to fibrations in \textbf{SympGr}. 
However, this is not the end of the story.

We have not used above the fact that the action $G\times \Sigma(F)\to\Sigma(F)$
is hamiltonian. By some standard results in symplectic geometry \cite{GLS,MaSa}, 
this implies that the fibered symplectic groupoid $\Sigma(M)\tto M$ carries a coupling 
2-form. Moreover, this closed 2-form is multiplicative, so that $\Sigma(M)$ is, in
fact, a presymplectic groupoid. This is also related with the problem of existence of
a coupling (see \cite{Vor}) for the original Poisson fibration. All this and related problems 
on Poisson fibrations is the subject of ongoing work (see \cite{BrFer}).

\subsection{Further constructions}

The symplectization functor can be (and should be!) applied to many other constructions
in Poisson geometry. 

For example, there exists a connected sum construction in Poisson
geometry (see \cite{IM}) which, under some conditions, out of two Poisson manifolds
$M_1$ and $M_2$ yields a new Poisson manifold $M_1\# M_2$. Some result
of the sort:
\[ \Sigma(M_1\# M_2)=\Sigma(M_1)\#\Sigma(M_2),\]
should be true (here, on the right-hand side one has a symplectic connected sum). This
kind of result should be relevant in the study of Poisson manifolds of compact type 
(\cite{CrFeMa}), namely in the problem of finding a decomposition into simpler pieces 
for such Poisson manifolds.

Another example, is in the theory of Poisson-Nijenhuis manifolds where the application of 
the $\Sigma$ functor leads to a symplectic-Nijenhuis groupoid (\cite{Cr,StXu}), and this
should be relevant in the study of the canonical integrable hierarchies associated with
a PN-manifold.

I believe that the symplectization functor $\Sigma$, which we have just started understanding, 
will play an important role in many other problems in Poisson geometry (\footnote{In the words of a
famous mathematician, ``This is no joke!''.}).

\section*{Acknowledgments}
I would like to thank the organizers of WGP for providing me with 
an excuse to write this article. I also thank 
L{\'\i}gia Abrunheiro, Rogier Bos, Raquel Caseiro, Jesus Clemente,  
Eva Miranda, Miguel Olmos and Patr{\'\i}cia Santos for many useful discussions
during the conference in Tenerife. This work has been partially supported by 
FCT/POCTI/FEDER and by grants POCI/MAT/55958/2004 and POCI/MAT/57888/2004.

\end{document}